\documentclass[envcountsame]{llncs}
\usepackage{amsmath}
\usepackage{amsfonts}
\usepackage{url}
\usepackage{color}
\usepackage{todonotes}
\usepackage[utf8]{inputenc}
\usepackage{hyperref}
\usepackage{graphicx}
\usepackage{algorithmic}
\usepackage{algorithm}
\usepackage{mathtools}

\usepackage{amsthm}
\DeclarePairedDelimiter{\ceil}{\lceil}{\rceil}
\DeclarePairedDelimiter\floor{\lfloor}{\rfloor}

\begin{document}

\title{Improved approximation of layout problems on random graphs}
\author{
Kevin K.~H.~Cheung\inst{1} \and
Patrick Girardet\inst{2}
}

\institute{
School of Mathematics and Statistics,
Carleton University,
Ottawa, ON, Canada.
\email{kevin.cheung@carleton.ca}
\and
Department of Mathematics, 
Rice University,
Houston, TX, USA.
\email{pdg2@rice.edu}
}

\maketitle

\begin{abstract}
Inspired by previous work of Diaz, Petit, Serna, and Trevisan
(Approximating layout problems on random graphs
Discrete Mathematics, 235, 2001, 245--253),
we show that several well-known graph layout problems are approximable to
within a factor arbitrarily close to 1 of the optimal with high probability for
random graphs drawn from an Erdös-Renyi distribution with appropriate sparsity
conditions.  Moreover, we show that the same results hold for the analogous
problems on directed acyclic graphs.
\end{abstract}

\section{Introduction}
Many well-known optimization problems on graphs fall into the category of
\textit{graph layout problems}. A layout of a graph on $n$ vertices consists of
a bijection between the vertices of the graph and the set $\{1,2,...,n\}$ of
natural numbers, which can be interpreted as arranging the vertices of the
graph in some order on a line. A graph layout problem then consists of
optimizing some objective function over the set of possible layouts of a graph.
There is an analogous notion of layout and layout problem for directed acyclic
graphs, wherein a layout of a directed acyclic graph is simply a topological
sort of it, so that the layout respects edge directions. The particular layout
problems we consider in this paper are the \textit{Minimum Cut Linear
Arrangement} (also known as \textit{Cutwidth}), \textit{Vertex Separation},
\textit{Edge Bisection}, and \textit{Vertex Bisection} problems, along with the
analogous problems on directed acyclic graphs. These problems are well known to
find applications in VLSI design, job scheduling, parallel computing, graph
drawing, etc. We direct the interested reader to a survey \cite{survey} on the
topic.

Graph layout problems are often computationally difficult to
solve exactly. Both the undirected and directed vertex separation problems are
known to be NP-complete \cite{vsepnp}, as are the undirected \cite{mincutnp}
and directed \cite{dmincutnp} minimum cut linear arrangement problems, and the
vertex bisection \cite{vbisnp} and edge bisection (shown in \cite{ebisnp} as a
special case of \textit{minimum cut into bounded sets}) problems on undirected
graphs. We do not know of a reference which proves the NP-completeness of the
vertex and edge bisection problems on directed graphs, though we have no reason
to believe that they are not also NP-complete. Due to the practical
applications of the problems considered, many researchers have sought to find
approximation algorithms for these problems. It is common to analyze the
performance of algorithms on random instances as a proxy for their ``real''
performance, so that one might seek to analyze the approximability of layout
problems on random graphs. 
Diaz {\it et al.} \cite{layout_approx} show that for any of
the undirected layout problems considered above, for any constant $C > 2,$ then
for large enough random graphs with appropriate sparsity conditions that any
solution of the problem has cost within a factor $C$ of the optimal with high
probability. Hence, these problems can be trivially approximated to within any
factor of $C>2$ for large enough random graphs with high probability. 

In this paper, in addition to showing that the constant of approximation
can be improved to any $C > 1$ with slightly weaker sparsity and convergence
results, we show that the same result holds for the directed versions of the
problems which were not considered in \cite{layout_approx}.
Moreover, we only use the Hoeffding inequality for tail bounds of sums of 
independent and identically distributed (i.i.d) random
variables and some well-known asymptotic estimates to show these results,
thus avoiding the more technical ``mixing graph'' framework used 
in \cite{layout_approx}.  In summary, for large enough random graphs
with appropriate sparsity conditions, any solution of these layout problems
will have cost arbitrarily close to optimal with high probability.   

\section{Definitions}
We first recall some terminology in \cite{layout_approx}.
Given an undirected graph $G=(V,E)$ with $|V| = n,$
a \textit{linear arrangement} (or a \textit{layout}) of $G$
is an injective function $\phi: V\rightarrow
\{1,...,n\}$ 
The problems we consider all take the form
of optimizing some objective function over the set of linear arrangements of a
graph. For a linear arrangement $\phi$ of $G,$ we have the two sets
$$L(i,\phi, G) = \{u\in V | \phi(u)\le i\}, R(i,\phi, G) = \{u\in V | \phi(u) >
i\}$$ 
and the two measures 
\begin{equation} \label{blah}
\begin{split}
\theta(i, \phi, G) & = |\{uv\in E | u\in L(i,\phi, G) \wedge v\in R(i,\phi,
G)\}|, \\ \delta(i, \phi, G) & = |\{u\in L(i,\phi, G) | \exists v\in
R(i,\phi,G) \wedge uv\in E\}|.
\end{split}
\end{equation}

The problems we consider for undirected graphs are the following: 
\begin{itemize}
\item Cutwidth (CUTWIDTH): Given a graph $G=(V,E),$ compute MINCW($G$) = $\min_\phi \text{CW}(G, \phi),$ where $\text{CW}(G, \phi) = \max_{i=1}^n \theta(i, \phi, G).$ 
\item Vertex separation (VERTSEP): Given a graph $G=(V,E),$ compute MINVS($G$) = $\min_\phi \text{VS}(G, \phi),$ where $\text{VS}(G, \phi) = \max_{i=1}^n \delta(i, \phi, G).$
\item Edge bisection (EDGEBIS): Given a graph $G=(V,E),$ compute MINEB($G$) = $\min_\phi \text{EB}(G, \phi),$ where $\text{EB}(G, \phi) = \theta(\floor{\frac{n}{2}}, \phi, G).$
\item Vertex bisection (VERTBIS): Given a graph $G=(V,E),$ compute MINVB($G$) = $\min_\phi \text{VB}(G, \phi),$ where $\text{VB}(G, \phi) = \delta(\floor{\frac{n}{2}}, \phi, G).$ 
\end{itemize}

These problems are all known to be NP-hard.

We also consider analogous problems on directed graphs. Given a directed
acyclic graph $G=(V,E)$ with $|V| = n,$ a \textit{layout} or \textit{linear
arrangement} of $G$ is an injective function $\phi: V\rightarrow \{1,...,n\}$
such that if $(u,v)\in E$ is a directed edge from $u\in V$ to $v\in V$ then
$\phi(u) < \phi(v).$ Note that this is simply a topological sort of $G,$ which
exists as $G$ is directed acyclic. The definitions of $L(i,\phi, G), R(i,\phi,
G), \theta(i, \phi, G), \delta(i, \phi, G)$ remain unchanged when switching to
directed acyclic graphs, and we can define the following directed versions of
the problems above: \begin{itemize}
\item Directed cutwidth (DCUTWIDTH): Given a directed acyclic graph $G=(V,E),$
compute MINCW($G$) = $\min_\phi \text{CW}(G, \phi),$ where $\text{CW}(G, \phi)$
is defined as above. The only difference is that we optimize over linear
arrangements of directed acyclic graphs instead.  \item Directed vertex
separation (DVERTSEP): Given a directed acyclic graph $G=(V,E),$ compute
MINVS($G$) = $\min_\phi \text{VS}(G, \phi),$ where
$\text{VS}(G, \phi) = \max_{i=1}^n \delta(i, \phi, G).$
\item Directed edge bisection (DEDGEBIS): Given a directed acyclic graph $G=(V,E),$ compute MINEB($G$) = $\min_\phi \text{EB}(G, \phi),$ where $\text{EB}(G, \phi) = \theta(\floor{\frac{n}{2}}, \phi, G).$
\item Directed vertex bisection (DVERTBIS): Given a directed acyclic graph $G=(V,E),$ compute MINVB($G$) = $\min_\phi \text{VB}(G, \phi),$ where $\text{VB}(G, \phi) = \delta(\floor{\frac{n}{2}}, \phi, G).$ 
\end{itemize}

For each arrangement problem considered above, we also define the maximum cost
solution of the problem on a graph. For example, for CUTWIDTH, in contrast to
MINCW($G$), we define MAXCW($G$) = $\max_\phi \text{CW}(G, \phi)$, and
similarly for every other problem considered above. Moreover, we define the gap
of a problem on a given graph $G$ to be the ratio of the highest cost solution
to the lowest cost solution. For example, for CUTWIDTH, the gap is
$$\text{GAPCW}(G) = \frac{\text{MAXCW}(G)}{\text{MINCW}(G)},$$ and this
quantity is defined in the same way for every other arrangement problem
considered above.

Any discussion on random graphs requires a probability distribution on graphs
In this paper, we adopt a variant of the Erdös-Renyi probability distribution
\cite{er_model} for undirected graphs defined as follows:

\textbf{Definition:} For a positive integer $n$ and probability $0\le p \le 1,$
the Erdös-Renyi distribution $G(n, p)$ on the set of $n$-vertex graphs assigns
an $n$-vertex graph $G=(V,E)$ probability $p^m (1-p)^{\binom{n}{2}} - m,$ where
$|E| = m.$ That is, we sample $n$-vertex graphs by including each possible edge
with probability $p.$ 

We also need a probability distribution on directed acyclic graphs. We use a
variant of the Erdös-Renyi probability distribution \cite{dag_model} which
produces directed acyclic graphs, defined as follows:

\textbf{Definition:} For a positive integer $n$ and probability $0\le p \le 1,$ the distribution $D(n, p)$ on the set of $n$-vertex directed acyclic graphs first samples a random graph from $G(n,p)$ on the vertex set $\{1,...,n\}.$ Then, each edge $\{i,j\}$ in the sampled graph is directed from $i$ to $j$ if $i < j$ and vice versa. 

As the edges in the sampled directed graph always point from a lower numbered
vertex to a higher numbered vertex, it is clear that the sampled graph is
acyclic.

\section{Preliminary lemmas}
We first list some technical lemmas necessary for carrying out the
probabilistic analysis in our main theorems. 
\begin{lemma}\label{hoeff}
[\textbf{Hoeffding's inequality:}] Suppose that $X_1,...,X_n$ are independent identically distributed Bernoulli random variables with mean $p,$ and let $H(n) = \overline{X_1}+\overline{X_2}+...+\overline{X_n},$ where $\overline{X_i}$ is a sample of $X_i.$ Then for $\epsilon > 0,$
\begin{align*}
    P[H(n) \le (p-\epsilon)n] &\le \exp(-2\epsilon^2n)\\
    P[H(n) \ge (p+\epsilon)n] &\le \exp(-2\epsilon^2n).
\end{align*}
\end{lemma}
\begin{proof}
This is a special case of Theorem 1 in Hoeffding's original paper \cite{hoeffding} for Bernoulli random variables. 
\end{proof}

\begin{lemma}\label{binomasy}
If $k = o(n),$ then $$\log\binom{n}{k}=(1+o(1))k\log\frac{n}{k}$$
\end{lemma}
\begin{proof}[Proof of \ref{binomasy}]
Recall the well-known inequalities $$(\frac{n}{k})^k \le \binom{n}{k} \le
(\frac{ne}{k})^k,$$ which can be obtained via Stirling's approximation. Taking
logs, we conclude that $$k\log\frac{n}{k}\le \log\binom{n}{k}\le
k(\log\frac{n}{k} + 1).$$ If $k = o(n),$ then
$\log\frac{n}{k}\rightarrow\infty,$ so that by the above chain of inequalities
we have that $$\log\binom{n}{k} = (1+o(1))k\log\frac{n}{k},$$ as desired.
\end{proof}

\begin{lemma}\label{cenbinom}
$$\binom{n}{\floor{n/2}}=\Theta(\frac{2^n}{\sqrt{n}}).$$
\end{lemma}
\begin{proof}
Recall Stirling's approximation that $$n! \sim \sqrt{2\pi n}(\frac{n}{e})^n.$$
By definition, $$\binom{2n}{n} = \frac{(2n)!}{(n!)^2},$$ so that
$$\binom{2n}{n}\sim \frac{\sqrt{4\pi n}(\frac{2n}{e})^{2n}}{2\pi
n(\frac{n}{e})^{2n}} = \frac{4^n}{\sqrt{\pi n}}.$$ Replace $n$ with $n/2$ to
obtain the desired result.
\end{proof}

\begin{lemma}\label{qest}
Suppose $f(n) = \Omega(n^{-c}), g(n) = \Omega(n^d),$ where $0<c < d.$ Moreover, suppose $f(n) = o(1).$ Then $$\lim_{n\to\infty}(1-f(n))^{g(n)} = 0.$$
\end{lemma}
\begin{proof}
Taking logarithms, we find that $$\log(\lim_{n\to\infty}(1-f(n))^{g(n)}) =
\lim_{n\to\infty}g(n)\log(1-f(n)) \le
\lim_{n\to\infty}k_1n^d\log(1-k_2n^{-c})$$ for appropriate constants $k_1, k_2
> 0.$ But then by L'Hopital's rule, $$\lim_{n\to\infty}k_1n^d\log(1-k_2n^{-c})
= \lim_{n\to\infty}\frac{\log(1-k_2n^{-c})}{k_1n^{-d}} =
\lim_{n\to\infty}\frac{-ck_2n^{-c-1}}{(1-k_2n^{-c})dk_1n^{-d-1}}.$$ Since
$0<c<d,$ we find that $$\frac{ck_2n^{-c-1}}{(1-k_2n^{-c})dk_1n^{-d-1}} =
O(n^{d-c})\rightarrow\infty,$$ so that $$\log(\lim_{n\to\infty}(1-f(n))^{g(n)})
= \lim_{n\to\infty}\frac{-ck_2n^{-c-1}}{(1-k_2n^{-c})dk_1n^{-d-1}} = -\infty,$$
and hence $$\lim_{n\to\infty}(1-f(n))^{g(n)} = 0,$$ as desired. 
\end{proof}

\section{Main results}
For the theorems that follow, let $\{G_n\}_{n=1}^\infty$ be a sequence of graphs sampled from a $G(n,p_n)$ Erdös-Renyi distribution with for some sequence $\{p_n\}_{n=1}^\infty$ of edge probabilities. 
\begin{theorem} \label{gapcw}
Let $p_n$ satisfy $p_n = \Omega(n^{-c}), c<1/2$. Then for all $\epsilon > 0, \delta > 0$ there exists an $N$ such that for all $n\ge N,$ $$P[\emph{GAPCW}(G_n) < 1 + \delta] > 1-\epsilon.$$ 
\end{theorem}

To establish this theorem we need the following lemmas: 
\begin{lemma} \label{lem1}
Let $p_n$ satisfy $p_n = \Omega(n^{-c}), c<1/2$. Then for all $\epsilon > 0, \delta > 0$ there exists an $N$ such that for all $n\ge N,$ $$P[\emph{MINCW}(G_n)\ge \frac{n^2p_n(1-\delta)}{4}] > 1-\epsilon.$$
\end{lemma}
\begin{lemma} \label{lem2}
Let $p_n$ satisfy $p_n = \Omega(n^{-c}), c<1/2$. Then for all $\epsilon > 0, \delta > 0$ there exists an $N$ such that for all $n\ge N,$ $$P[\emph{MAXCW}(G_n)\le \frac{n^2p_n(1+\delta)}{4}] > 1-\epsilon.$$
\end{lemma}
\begin{proof}[Proof of Lemma \ref{lem1}]
Suppose we are given an Erdös-Renyi random graph $G_n = (V,E)$ with $n$ vertices and edge probability $p_n$. For a partition of $V$ into disjoint sets $A, V-A,$ let $c(A)$ denote the number of edges going between $A$ and $V-A.$ That is, 
$$c(A) = |\{uv\in E\;s.t.\;u\in A, v\in V-A\}|.$$ 
Given a linear arrangement $\phi: V\rightarrow [n]$ we write 
$$L_\phi = L(\floor*{n/2}, \phi, G) = \{u\in V\;s.t.\; \phi(u)\le \floor*{n/2}\},$$ 
$$R_\phi = R(\floor*{n/2}, \phi, G) = \{u\in V\;s.t.\; \phi(u) > \floor*{n/2}\}$$ to be the sets of vertices respectively in the left and right halves of the linear arrangement. It is clear then that $\theta(G_n, \phi) \ge c(L_\phi),$ so that MINCW($G_n$) $\ge \min_\phi c(L_\phi).$ This last expression is equivalent to the minimum of $c(S)$ over subsets $S$ of $V$ with cardinality $|S| = \floor{n/2}.$ 

Let $S$ be a random subset of $V$ of size $|S| = \floor{n/2}.$ As $G_n$ is an
Erdös-Renyi random graph with edge probability $p_n, c(S)$ is a binomial random
variable with mean $\mu = \floor{n/2}\ceil{n/2}p_n \ge \frac{(n^2 - 1)p_n}{4}$.
By applying Hoeffding's inequality \ref{hoeff} to $c(S),$ we have that
$$P[c(S)\le \mu(1-\frac{\alpha_n}{p_n})]\le \exp(-(n^2 - 1)\alpha_n^2/2)$$ for
$\alpha_n > 0.$ As $p_n = \Omega(n^{-c})$ with $c>1/2,$ we can choose
$\alpha_n$ to get the desired convergence: set $\alpha_n = n^{-l}$ where $l$
satisfies $c<l<1/2,$ so that $\alpha_n = o(p_n)$ and  $\alpha_n^2 =
\Omega(n^{-1+s})$ for some $s>0.$ We thus have that \begin{align*}
    P[c(S) &\le \mu(1-\frac{\alpha_n}{p_n})]\le \exp(-(n^2 - 1)\alpha_n^2/2) = \exp(-(n^2 - 1)\Omega(n^{-1+s}))\\
           &= \exp(-\Omega(n^{1+s}))
\end{align*}     
for some $s>0.$ Note also that $\dfrac{\alpha_n}{p_n}=O(n^{c-l})=o(1),$ so that
the above can be written as $$P[c(S) \le \mu(1-O(n^{c-l}))]\le
\exp(-\Omega(n^{1+s})).$$ 
Note trivially that $\dbinom{n}{\floor{n/2}}\le 2^n.$ Thus, applying the union
bound to this inequality over all $S\subset V$ with $|S|=\floor{n/2},$ we find
that the probability that \textbf{any} such $S$ has cut
$c(S)\le\mu(1-O(n^{c-l}))$ is at most
$$\dbinom{n}{\floor{n/2}}P[c(S)\le\mu(1-O(n^{c-l}))]\le
2^n\exp(-\Omega(n^{1+s}))\rightarrow 0$$ as $n\rightarrow\infty.$ Thus, by our
previous discussion $$\text{MINCW}(G_n)\ge\min_\phi c(L_\phi)\ge
\mu(1-O(n^{c-l}))\ge \frac{(n^2-1)p_n(1-O(n^{c-l}))}{4}$$ with probability at
least $1-2^n\exp(-\Omega(n^{1+s}))$. Then for any $\delta>0,$ for large enough
$n$ we have that $$\frac{(n^2-1)p_n(1-O(n^{c-l}))}{4} \ge
\frac{n^2p_n(1-\delta)}{4},$$ and for any $\epsilon > 0,$ for large enough $n$
we have that $$2^n\exp(-\Omega(n^{1+s})) < \epsilon,$$ so that
$$\text{MINCW}(G_n)\ge \frac{n^2p_n(1-\delta)}{4}$$ with probability at least
$1-\epsilon$ for large enough $n,$ as desired. 
     
\end{proof}
\begin{proof}[Proof of Lemma \ref{lem2}]
Let $\phi$ be a linear arrangement of $G_n$ such that $\text{CW}(G_n,\phi)
= k.$ It must then be the case that there exists some $S\subset V$ (no
cardinality constraint this time) such that $c(S)=k,$ since by definition there
must exist some $i$ such that $\text{CW}(G_n,\phi) = \theta(i,\phi,G_n),$ and
we just take $S=L(i,\phi,G), V-S=R(i,\phi,G)$ as defined previously. Hence, an
upper bound with high probability on $c(S)$ for all $S\subset V, |S|\le
\floor{n/2}$ also gives an upper bound with high probability on
$\text{MAXCW}(G_n).$ We thus follow a similar proof as that for lemma 1. Since
$f(k)=k(n-k)$ is maximal for $k=n/2,$ if $S$ is a random subset of $V$ with
$|S|\le\floor{n/2},$ then for any $c$ we have that $P[c(S)\ge c]\le
P[c(X)\ge c]$ if $X$ is a random subset of $V$ with $|X|=\floor{n/2}.$ Thus, by
the Hoeffding inequality, if $S$ is a random subset of $V$ with
$|S|\le\floor{n/2},$ then $$P[S\ge\mu(1+\frac{\alpha_n}{p_n})]\le \exp(-(n^2 -
1)\alpha_n^2/2)$$ for
$\alpha_n > 0$ where $\mu = \floor{n/2}\ceil{n/2}p_n$ as above. As $p_n =
\Omega(n^{-c})$ with $c>1/2,$ we again set $\alpha_n = n^{-l}$ where $l$
satisfies $c<l<1/2,$ so that $\alpha_n = o(p_n)$ and  $\alpha_n^2 =
\Omega(n^{-1+s})$ for some $s>0.$ As before, we have that
$$P[c(S) \ge \mu(1+\frac{\alpha_n}{p_n})] \le \exp(-\Omega(n^{1+s}))$$     
for some $s>0,$ which we can again write as $$P[c(S) \ge \mu(1+O(n^{c-l}))]\le
\exp(-\Omega(n^{1+s})).$$ 
Hence, by the union bound, the probability that any $S\subset V,
|S|\le\floor{n/2}$ has $c(S)\ge\mu(1+O(n^{c-l}))$ is at most
$$2^nP[c(S)\ge\mu(1+O(n^{c-l}))]\le 2^n\exp(-\Omega(n^{1+s}))\rightarrow 0$$ as
$n\rightarrow\infty.$ Thus, $\text{MAXCW}(G_n)\le
\dfrac{n^2p_n(1+O(n^{c-l}))}{4}$ with probability at least
$1-2^n\exp(-\Omega(n^{1+s})).$  Thus, as before, for any $\delta>0,$ for large
enough $n$ we have that $$\frac{n^2p_n(1+O(n^{c-l}))}{4}\le
\frac{n^2p_n(1+\delta)}{4},$$ and for any $\epsilon >0,$ for large enough $n$
we have that $$2^n\exp(-\Omega(n^{1+s}))<\epsilon,$$  so that
$$\text{MAXCW}(G_n)\le \frac{n^2p_n(1+\delta)}{4}$$ with probability at least
$1-\epsilon$ for large enough $n,$ as desired.
\end{proof}

With these two lemmata, the main theorem for the cutwidth gap is easy to
establish. 
\begin{proof}[Proof of Theorem \ref{gapcw}]
As in the statement of the theorem, let $p_n$ satisfy $p_n = \Omega(n^{-c})$
for some $c<1/2,$ and let $\delta, \epsilon > 0$ be given. Since $$\lim_{x\to
0}\frac{1+x}{1-x} = 1,$$ let $\delta'$ satisfy $\dfrac{1+\delta'}{1-\delta'} <
1 + \delta.$ By lemma \ref{lem1}, pick $N_1$ so that for $n\ge N_1,$
$$P[\text{MINCW}(G_n)\ge \frac{n^2p_n(1-\delta')}{4}] > 1-\epsilon/2,$$ and
similarly by lemma \ref{lem2} pick $N_2$ so that for $n\ge N_2,$
$$P[\text{MAXCW}(G_n)\le\frac{n^2p_n(1+\delta')}{4}] > 1-\epsilon/2.$$ Hence,
if $N=\max\{N_1,N_2\},$ then for $n\ge N$ we have that \begin{align*}
    1 - \epsilon &= 1-2(\epsilon/2) \\
     &\le P[(\text{MINCW}(G_n)\ge \frac{n^2p_n(1-\delta')}{4})\wedge(\text{MAXCW}(G_n)\le\frac{n^2p_n(1+\delta')}{4})].
\end{align*}
Since if $$\text{MINCW}(G_n)\ge \frac{n^2p_n(1-\delta')}{4}$$ and
$$\text{MAXCW}(G_n)\le\frac{n^2p_n(1+\delta')}{4}$$ we must have that
$$\text{GAPCW}(G_n) = \frac{\text{MAXCW}(G_n)}{\text{MINCW}(G_n)} \le
\frac{1+\delta'}{1-\delta'} < 1 + \delta,$$ and hence that if $n\ge N,$ then
$$1-\epsilon < P[\text{GAPCW}(G_n) < \frac{1+\delta'}{1-\delta'}] \le
P[\text{GAPCW}(G_n) < 1+\delta],$$ as desired.
\end{proof}

We then obtain as a corollary the analogous gap result for directed graphs. We
should of course expect such a result, since a random directed acyclic graph
drawn from the $D(n,p)$ model is essentially a $G(n,p)$ graph but with fewer
valid linear arrangements, allowing for less variance in $\text{CW}(G,\phi).$
\begin{theorem} \label{gapdcw}
Let $p_n$ satisfy $p_n = \Omega(n^{-c}), c<1/2,$ and let $\{G_n\}$ be a sequence of directed acyclic graphs sampled from a $D(n,p_n)$ distribution. Then for all $\epsilon > 0,\delta>0$ there exists an $N$ such that for all $n\ge N,$ $$P[\emph{GAPCW}(G_n)<1+\delta] > 1-\epsilon.$$
\end{theorem}
\begin{proof}
Let $G = (V,E)$ be a directed acyclic graph. For a partition of $V$ into
disjoint sets $A, V-A,$ let $c(A)$ denote the number of edges going from $A$ to
$V-A.$ That is, $$c(A) = |\{(u,v)\in E \; s.t. \; u\in A, v\in V-A\}|.$$ As in
the proof of lemma \ref{lem1}, we can define $L_\phi, R_\phi$ for a linear
arrangement $\phi:V\rightarrow [n],$ and it is again clear that
$\text{MINCW}(G_n)\ge\min_\phi c(L_\phi).$ Since not all permutations of $V$
are necessarily valid linear arrangements for the directed acyclic graph $G,$
this last expression is a minimization over at most $\dbinom{n}{\floor{n/2}}$
subsets $S$ of $V$ with cardinality $|S|=\floor{n/2},$ but in general fewer.
Now if $S$ is a subset of $V$ which arises as $L_\phi$ for some linear
arrangement $\phi,$ then by the definition of $D(n,p)$ we have that $c(S)$ is a
binomial random variable with mean $\mu = \floor{n/2}\ceil{n/2}p_n\ge
\frac{(n^2-1)p_n}{4}.$ We can then apply Hoeffding's inequality and emulate the
analysis in the proof of lemma \ref{lem1} to once again conclude that for large
enough $n$ that $$\text{MINCW}(G_n)\ge\frac{n^2p_n(1-\delta)}{4}$$ with
probability at least $1-\epsilon$ for any $\delta, \epsilon > 0.$ The proof
only differs in that we note that there are at most (rather than exactly)
$\dbinom{n}{\floor{n/2}}$ subsets of $V$ arising as $L_\phi$ for some linear
arrangement $\phi,$ so that our application of the union bound may be wasteful
in general, which only strengthens our result.

We can of course also essentially copy the proof of lemma \ref{lem2},
restricting our attention to subsets $S\subset V$ which are of the form
$L(i,\phi,G)$ for some $i$ and valid linear arrangement $\phi.$ For subsets $S$
of the form $L_\phi$ for some linear arrangement $\phi,$ we can obtain the same
upper bound $$P[c(S)\ge\mu(1+\frac{\alpha_n}{p_n})]\le\exp(-\Omega(n^{1+s})),$$
and then apply the union bound in the same way (which may be wasteful as not
all subsets of $V$ are necessarily of the form $L(i,\phi,G).$ The analysis then
carries through verbatim, so that
$$\text{MAXCW}(G_n)\le\frac{n^2p_n(1+\delta)}{4}$$ with probability at least
$1-\epsilon$ for large $n$ as before.

Combining the above two probabilistic inequalities precisely as in the proof of
theorem \ref{gapcw}, we have that for any $\delta,\epsilon > 0,$
$$P[\text{GAPCW}(G_n)<1+\delta]>1-\epsilon$$ for large enough $n,$ giving the
desired analogue of theorem \ref{gapcw} for random directed acyclic graphs, so
that DCUTWIDTH exhibits the same convergence behavior as CUTWIDTH. 
\end{proof}

Moreover, we can further emulate the proof of theorem \ref{gapcw} to obtain
analogous convergence results for EDGEBIS, DEDGEBIS.

\begin{theorem}
\label{gapeb}
Let $p_n$ satisfy $p_n=\Omega(n^{-c}), c<1/2,$ and let be $\{G_n\}$ be a sequence of random graphs sampled from $G(n,p_n).$ Then for all $\epsilon > 0, \delta > 0$ there exists an $N$ such that for all $n\ge N,$ 
$$P[\emph{GAPEB}(G_n)<1+\delta]>1-\epsilon.$$
\end{theorem}
\begin{proof}
Note that $\text{MINEB}(G_n)$ is simply the minimization of $c(S)$ over subsets
$S\subset V$ with $|S|=\floor{n/2},$ and $\text{MAXEB}(G_n)$ the corresponding
maximization. Hence, the proof of lemma \ref{lem1} carries through without
modification to yield that $$P[\text{MINEB}(G_n)\ge
\frac{n^2p_n(1-\delta)}{4}]>1-\epsilon$$ for any $\epsilon,\delta > 0$ for
large $n.$ Similarly, the proof of lemma \ref{lem2} carries through to yield
that $$P[\text{MINEB}(G_n)\ge \frac{n^2p_n(1+\delta)}{4}]>1-\epsilon$$ for any
$\epsilon,\delta > 0$ for large $n,$ with only a slightly more wasteful union
bound application since we need only consider subsets $S\subset V$ of size
$|S|=\floor{n/2}.$ Combining these two results in the same manner as in the
proof of theorem \ref{gapcw} gives the desired result.
\end{proof}

\begin{theorem} \label{gapdeb}
Let $p_n$ satisfy $p_n=\Omega(n^{-c}), c<1/2,$ and let be $\{G_n\}$ be a sequence of random directed acyclic graphs sampled from $D(n,p_n).$ Then for all $\epsilon > 0, \delta > 0$ there exists an $N$ such that for all $n\ge N,$ 
$$P[\emph{GAPEB}(G_n)<1+\delta]>1-\epsilon.$$
\end{theorem}
\begin{proof}
Simply carry through the proof of theorem \ref{gapeb}, observing that for a
directed acyclic graph $G,$ we only need to consider subsets $S\subset V$ of
the form $L_\phi$ for a valid linear arrangement $\phi,$ rather than all
$S\subset V$ with $|S|=\floor{n/2}.$ As in the proof of theorem \ref{gapdcw},
this will give slightly wasteful union bound applications in the proofs of the
probabilistic lower and upper bounds for $\text{MINEB}(G_n),
\text{MAXEB}(G_n),$ but the same bounds will hold which can then be combined
verbatim as in the proof of theorem \ref{gapcw} to give the desired result.
\end{proof}

Thus, all of our graph arrangement problems involving overhead edges, that is,
CUTWIDTH, DCUTWIDTH, EDGEBIS, DEDGEBIS, all exhibit the same convergence
behavior for large random graphs. We now turn our attention to proving
analogous convergence results for our vertex arrangement problems, that is,
VERTSEP, DVERTSEP, VERTBIS, DVERTBIS.  \begin{theorem}\label{gapvs}
Let $p_n$ satisfy $p_n=\Omega(n^{-c}), c<1,$ and let $\{G_n\}$ be a sequence of random graphs sampled from $G(n,p_n).$ Then for all $\epsilon > 0, \delta > 0$ there exists an $N$ such that for all $n\ge N,$ 
$$P[\emph{GAPVS}(G_n)<1+\delta]>1-\epsilon.$$
\end{theorem}

The following lemma essentially proves the theorem. 

\begin{lemma}\label{lem3}
Let $p_n$ satisfy $p_n=\Omega(n^{-c}), c<1,$ and let $\{G_n\}$ be a sequence of random graphs sampled from $G(n,p_n).$ Then for all $\epsilon > 0, \delta > 0$ there exists an $N$ such that for all $n\ge N,$ $$P[\emph{MINVS}(G_n)\ge (1-\delta)n - 1] > 1-\epsilon.$$
\end{lemma}
\begin{proof}
Set $\delta_n = n^{-l}, \epsilon_n = n^{-s},$ where $0<s<1-c, 0<l<\frac{s}{2}.$
Given a linear arrangement $\phi$ of $G_n,$ set $S_{\phi,\epsilon_n} =
L(\floor{(1-\epsilon_n)n}, \phi, G_n).$ As there are $\ceil{\epsilon_n n}$
vertices in $V-S_{\phi,\epsilon_n},$ the odds of a given vertex $v\in
S_{\phi,\epsilon_n}$ not being connected to any vertex in
$V-S_{\phi,\epsilon_n}$ are $(1-p_n)^{\ceil{\epsilon_n n}}.$ Hence, the number
of vertices in $S_{\phi,\epsilon_n}$ which are connected to elements of
$V-S_{\phi,\epsilon_n}$ is a binomial random variable $X$ on
$m=\floor{(1-\epsilon_n)n}$ trials and probability
$q=1-(1-p_n)^{\ceil{\epsilon_n n}}.$ By Hoeffding's inequality, we have that
$$P[X\le(q-\delta_n)m]\le \exp(-2\delta_n^2m).$$ Any linear arrangement of
$G_n$ gives rise to such a set $S_{\phi,\epsilon_n},$ of which there are
$\dbinom{n}{\ceil{\epsilon_n n}}.$ Hence, by the union bound, the probability
that any linear arrangement of $G_n$ has a vertex separation of less than
$(q-\delta_n)m$ is at most 
$$\binom{n}{\ceil{\epsilon_n
n}}P[X\le(q-\delta_n)m]\le\binom{n}{\ceil{\epsilon_n n}}\exp(-2\delta_n^2 m).$$
By lemma \ref{binomasy}, we find that the log of the probability that any
linear arrangement of $G_n$ has vertex separation less than $(q-\delta_n)m$ is
at most $$\log\binom{n}{\ceil{\epsilon_n n}} + \log\exp(-2\delta_n^2 m) =
(1+o(1))\ceil{\epsilon_n n}\log{\frac{1}{\epsilon_n}}-2\delta_n^2m,$$ since
$\ceil{\epsilon_n n}=\ceil{n^{1-s}}=o(n).$ Substituting in our choices of
$\delta, \epsilon,$ we find that the log of the relevant probability is at most
$$(s+o(1))n^{1-s}\log(n)-2n^{1-2l}(1-n^{-s})\rightarrow -\infty$$ as
$n\rightarrow\infty,$ since $1-2l>1-s$ by the condition that $l<\frac{s}{2}.$
Thus, for large $n$ we have that $P[\text{MINVS}(G_n)\le(q-\delta_n)m] <
\epsilon.$ Moreover, if we expand out $q,\delta_n,m,$ we find that
$$(q-\delta_n)m=(1-(1-p_n)^{\ceil{n^{1-s}}}-n^{-l})\floor{(1-n^{-s})n}.$$ We
have by lemma \ref{qest} that $(1-(1-p_n)^{\ceil{n^{1-s}}} = o(1),$ so that 
$$(q-\delta_n)m=(1-(1-p_n)^{\ceil{n^{1-s}}}-n^{-l})\floor{(1-n^{-s})n} = (1-o(1))n - 1.$$ Hence, by combining these two results we find that for large $n$ we have that $$P[\text{MINVS}(G_n)\ge (1-\delta)n - 1] > 1-\epsilon,$$ as desired. 
\end{proof}

\begin{proof}[Proof of Theorem \ref{gapvs}]
Observe trivially that $\text{MAXVS}(G_n)\le n-1$ since $G_n$ only has $n-1$ vertices. Let $\epsilon, \delta>0$ be given. Set $\delta'$ so that for large $n, \dfrac{n-1}{(1-\delta')n-1}< 1+\delta.$ If we invoke lemma \ref{lem3} with $\delta', \epsilon,$ we have an $N$ such that for $n\ge N,$  $$P[\emph{MINVS}(G_n)\ge (1-\delta)n - 1] > 1-\epsilon.$$ Thus, 
$$\frac{\text{MAXVS}(G_n)}{\text{MINVS}(G_n)}\le\frac{n-1}{(1-\delta')n-1}<1+\delta$$ with probability greater than $1-\epsilon,$ so that $$P[\text{GAPVS}(G_n)<1+\delta]>1-\epsilon$$ as desired. 
\end{proof}

The analogous convergence results for DVERTSEP, VERTBIS, DVERTBIS
follow similarly.

\begin{theorem}\label{gapdvs}
Let $p_n$ satisfy $p_n=\Omega(n^{-c}), c<1,$ and let $\{G_n\}$ be a sequence of
random directed acyclic graphs sampled from $D(n,p_n).$ Then for all $\epsilon
> 0, \delta > 0$ there exists an $N$ such that for all $n\ge N,$
$$P[\emph{GAPVS}(G_n)<1+\delta]>1-\epsilon.$$
\end{theorem}
\begin{proof}
As in the proof of \ref{gapvs}, it suffices to prove that the conclusion of
lemma \ref{lem3} holds if $G_n$ is instead sampled from $D(n,p_n).$ The proof
of \ref{lem3} carries through verbatim with the one caveat that not all subsets
$S\subset V$ of size $|S|=\floor{(1-\epsilon_n)n}$ arise as
$S_{\phi,\epsilon_n}$ for some linear arrangement $\phi,$ so that our
application of the union bound is slightly wasteful. Otherwise, the proof goes
through, so that combining this lemma with the fact that $\text{MAXVS}(G_n)\le
n-1$ gives the desired result for random directed acyclic graphs.  \end{proof}

The analogous results for VERTBIS, DVERTBIS don't follow quite so immediately
as we have seen before but are still true.

\begin{theorem}\label{gapvb}
Let $p_n$ satisfy $p_n=\Omega(n^{-c}), c<1,$ and let $\{G_n\}$ be a sequence of random graphs sampled from $G(n,p_n).$ Then for all $\epsilon > 0, \delta > 0$ there exists an $N$ such that for all $n\ge N,$ 
$$P[\emph{GAPVB}(G_n)<1+\delta]>1-\epsilon.$$
\end{theorem}
\begin{proof}
Given a linear arrangement $\phi$ of $G_n,$ we have that $\text{VB}(G_n,\phi)$
is just the number of vertices in $L_\phi$ connected to vertices in $R_\phi.$
The number of vertices in $L_\phi$ connected to vertices in $R_\phi$ is a
binomial random variable $X$ on $m=\floor{n/2}$ trials with probability $q =
1-(1-p_n)^{\floor{n/2}}.$ By Hoeffding's inequality, we have that $$P[X\le
(q-\delta')m] = P[X\le (q-\delta)\frac{n}{2}] \le \exp(-\delta^2n).$$ Any
linear arrangement of $G_n$ gives rise to such a set $L_\phi,$ of which there
are $\dbinom{n}{\floor{n/2}}.$ Hence, by the union bound, the probability that
any linear arrangement of $G_n$ has a vertex separation of less than
$(q-\delta')m$ is at most
$$\binom{n}{\floor{n/2}}P[X\le(q-\delta')m]\le\binom{n}{\floor{n/2}}\exp(-\delta'^2n)=\Theta(\frac{2^n}{\sqrt{n}})\exp(-\delta'^2n)$$
by the estimate \ref{cenbinom} for the central binomial coefficient. For large
$n$ the $\sqrt{n}$ in the denominator of this last expression means that
$$\binom{n}{\floor{n/2}}P[X\le(q-\delta')m]\le\Theta(\frac{2^n}{\sqrt{n}})\exp(-\delta'^2n)<\epsilon.$$
Now since $(1-p_n)^{\floor{n/2}}\rightarrow 0$ since $p_n = \Omega(n^{-c}),$ we
have that for large $n$ and appropriately small choice of $\delta'$ that the
probability that any $L_\phi$ is connected to fewer than
$(1-\delta)\frac{n}{2}$ vertices in $R_\phi$ is at most $\epsilon,$ so that
$$P[\text{MINVB}(G_n)\ge (1-\delta)\frac{n}{2}] > 1-\epsilon.$$ Since
$\text{MAXVB}(G_n)\le \frac{n}{2}$ by definition, this gives the desired
result.
\end{proof}

\begin{theorem}\label{gapdvb}
Let $p_n$ satisfy $p_n=\Omega(n^{-c}), c<1,$ and let $\{G_n\}$ be a sequence of random directed acyclic graphs sampled from $D(n,p_n).$ Then for all $\epsilon > 0, \delta > 0$ there exists an $N$ such that for all $n\ge N,$ 
$$P[\emph{GAPVB}(G_n)<1+\delta]>1-\epsilon.$$
\end{theorem}
\begin{proof}
Simply repeat the proof of theorem \ref{gapvb}, and note that since not every subset $S\subset V$ with $|S|=n/2$ may necessarily arise as $L_\phi$ for some linear arrangement $\phi,$ the union bound might be wasteful. 
\end{proof}

Hence, we have shown that every problem we initially defined has
the desired gap convergence under appropriate conditions. 

\section{Concluding remarks}
In this paper, we have shown that many layout problems of interest can be
approximated arbitrarily close to the optimal with high probability for large
random graphs with appropriate sparsity conditions. This improves upon the
previous best results which demonstrated that some of these problems could be
approximated arbitrarily close to a factor of 2 times the optimal. We note that
there is still room for improvement with our results. The previous factor of 2
approximations in \cite{layout_approx} held for edge probabilities $p_n =
\Omega(n^{-1}),$ whereas our results for the layout problems on edges (i.e.
minimum cut linear arrangement, edge bisection, and directed versions) only
hold for $p_n = \Omega(n^{-c})$ for $c<1/2.$ Hence, improving the sparsity
conditions under which these approximability results hold is an open problem.
Moreover, the results in \cite{layout_approx} don't experience the same
tradeoff between sparsity and speed of convergence that our results do, a
seeming consequence of the strength of their ``mixing graph'' framework, so
that improving the speed of convergence also remains an area of potential
improvement. Some of the key results about mixing graphs used in
\cite{layout_approx} call upon the Hoeffding inequality, which was our primary
probabilistic tool in this paper. Hence, it would be interesting to see whether
the techniques of this paper and the mixing graphs 
could be unified somehow to give our improved constant of approximation but
retain the better sparsity and convergence conditions of \cite{layout_approx}. 

\section*{Acknowledgement}
The authors would like to thank Emmanuel Ruiz and Ashkan Moatamed for
conversations on research involving graph layout problems.
Additionally, the research in this paper was made possible by the support
the Fields Institute through its 2017 Fields Undergraduate Summer Research
Program.

\end{document}